# Hodge operators and groups of isometries of diagonalizable symmetric bilinear forms in characteristic two

Linus Kramer [a], Markus J. Stroppel [b,*]

[a] *Mathematisches Institut, Fachbereich Mathematik und Informatik, Universität Münster, Einsteinstraße 62, 48149 Münster, Germany*
[b] *LExMath, Fakultät 8, Universität Stuttgart, 70550 Stuttgart, Germany*

A R T I C L E   I N F O

A B S T R A C T

We study groups of isometries of non-alternating symmetric bilinear forms on vector spaces of characteristic two, and actions of these groups on exterior powers of the space, viewed as modules over algebras generated by Hodge operators.

* Corresponding author.
*E-mail addresses:* linus.kramer@uni-muenster.de (L. Kramer), stroppel@mathematik.uni-stuttgart.de (M.J. Stroppel).

## Introduction

In [10], we have introduced Hodge operators using diagonalizable $\sigma$-hermitian forms on vector spaces over a field $\boldsymbol{F}$ of arbitrary characteristic; here $\sigma$ is an automorphism of $\boldsymbol{F}$ with $\sigma^2 = \mathrm{id}$, we allow $\sigma = \mathrm{id}$ (in fact, we concentrate on that case in the present paper). If the dimension of that space is an even number $n = 2\ell$ then these operators help to understand exceptional homomorphisms between groups of semi-similitudes; these homomorphisms can be interpreted as representations of the group of semi-similitudes of the given form on the $\ell$-th exterior power of the space, where the latter is turned into a module over an algebra $\boldsymbol{K}_\ell$ generated by the Hodge operator (see [10, Sect. 2]). That algebra turns out to be a composition algebra if $\sigma \neq \mathrm{id}$ or if $\operatorname{char} \boldsymbol{F} \neq 2$, but it will be inseparable if $\sigma = \mathrm{id}$ and $\operatorname{char} \boldsymbol{F} = 2$: in that case, we obtain $\boldsymbol{K}_\ell \cong \boldsymbol{F}[X]/(X^2 - \delta)$, for some $\delta \in \boldsymbol{F}$.

In the case where $\sigma = \mathrm{id}$ and the characteristic is two, the forms in question are not the ones that lead to classical groups: we then use a bilinear form that is symmetric but not alternating (see 1.2 below); the definitions of classical groups (i.e., symplectic, unitary, or orthogonal groups) in characteristic two employ non-degenerate forms that are either alternating, or $\sigma$-hermitian with $\sigma \neq \mathrm{id}$, or quadratic forms.

Therefore, the inseparable case is treated in a cursory way only in [10]. However, it leads to phenomena that appear to be interesting, if only as a marked contrast to the results in [10]: for instance, the orthogonal group with respect to the bilinear form may be a simple group (see 3.2.a below) or abelian (see 3.2.b), or it may act on a submodule of the exterior power with a rather large nilpotent kernel, and inducing a simple group on that submodule (see 3.1). We treat this inseparable case in a more detailed manner in the present notes, with a focus on $\ell = 2$ (and $n = 4$) because interesting phenomena are already apparent in this dimension (and the Klein quadric provides some extra geometric intuition).

## 1. Symmetric bilinear and diagonal quadratic forms in characteristic two

We recall basic facts about forms over fields of characteristic two, and fix some notation.

**1.1 Notation.** Let $\boldsymbol{F}$ be any field of characteristic 2, let $V$ be a vector space of dimension $n$ over $\boldsymbol{F}$, and let $h\colon V \times V \to \boldsymbol{F}$ be a non-degenerate symmetric bilinear form. Moreover, assume that there exists an orthogonal basis $v_1, \dots, v_n$ with respect to $h$ (see 1.2 below).

As $\dim V$ is assumed to be finite, our assumption that $h$ be not degenerate is equivalent to the fact that (by a slight abuse of notation) we may view $h$ as a linear isomorphism onto the dual space $V^\vee$,

$$h\colon V \to V^\vee : v \mapsto h(v, -)$$

see e.g. [9, Ch. I, §2]. Consider now the exterior algebra $\bigwedge V$, cf. [8, VI 9]. We note that $\bigwedge$ is a functor on vector spaces and (semi)linear maps, cp. [10, 1.6]. Moreover, there is a natural isomorphism $(\bigwedge V)^\vee \cong \bigwedge(V^\vee)$, so we may write unambiguously $\bigwedge V^\vee$. Explicitly, we have $\langle f_1 \wedge \cdots \wedge f_\ell, w_1 \wedge \cdots \wedge w_\ell\rangle = \det(\langle f_i, w_k\rangle)$ for $f_i \in V^\vee$ and $w_k \in V$, see [11, I.5.6] or[1] [2, § 8, Thme. 1, p. 102]. Applying the functor $\bigwedge$ to $h\colon V \to V^\vee$, we obtain $\bigwedge h\colon \bigwedge V \to \bigwedge V^\vee$; we interpret this as a bilinear form $\bigwedge h$ on the exterior algebra $\bigwedge V$. Using the explicit formula above, we find

$$\bigwedge h(v_1 \wedge \cdots \wedge v_\ell, w_1 \wedge \cdots \wedge w_\ell) = \bigwedge^\ell h(v_1 \wedge \cdots \wedge v_\ell, w_1 \wedge \cdots \wedge w_\ell) = \det(h(v_i, w_j)).$$

In particular, the form $\bigwedge^\ell h$ is symmetric because transposition does not change the determinant.

The quadratic form $q \coloneqq q_h\colon V \to \boldsymbol{F}\colon v \mapsto h(v,v)$ is a semilinear map, its companion is the Frobenius endomorphism $\varphi\colon \boldsymbol{F} \to \boldsymbol{F}\colon x \mapsto x^2$, considered as an isomorphism from $\boldsymbol{F}$ onto the field $\boldsymbol{F}^\square \coloneqq \{x^2 \mid x \in \boldsymbol{F}\}$ of squares. The kernel of $q$ consists of all vectors that are isotropic with respect to $h$. In particular, it contains $V^\perp$. For each vector space complement $W$ for that kernel in $V$, the restriction of $q$ to $W$ is an injective map, and $\dim_{\boldsymbol{F}} W = \dim_{\boldsymbol{F}^\square} q(V)$.

Let $\Gamma\mathrm{L}(V)$ be the group of semilinear bijections of $V$ onto itself. For $\gamma \in \Gamma\mathrm{L}(V)$, let $\alpha_\gamma \in \mathrm{Aut}(\boldsymbol{F})$ denote the companion of $\gamma$. In the spirit of Dieudonné's notation, we consider the group

$$\Gamma\mathrm{O}(V,h) \coloneqq \left\{\gamma \in \Gamma\mathrm{L}(V) \,\middle|\, \exists r_\gamma \in \boldsymbol{F}^\times\ \forall v,w \in V\colon h(\gamma(v),\gamma(w)) = r_\gamma\, \alpha_\gamma(h(v,w))\right\}$$

of *semi-similitudes* of the form $h$, the subgroup $\mathrm{GO}(V,h) \coloneqq \Gamma\mathrm{O}(V,h) \cap \mathrm{GL}(V)$ of (linear) similitudes, and the group $\mathrm{O}(V,h) \coloneqq \left\{\gamma \in \mathrm{GL}(V) \,\middle|\, \forall v,w \in V\colon h(\gamma(v),\gamma(w)) = h(v,w)\right\}$ of isometries. Note that Dieudonné would refer to the latter group as a unitary group, he reserved the term "orthogonal" for groups of isometries of quadratic forms.

For a quadratic form $q\colon V \to \boldsymbol{F}$ we have

$$\Gamma\mathrm{O}(V,q) \coloneqq \left\{\gamma \in \Gamma\mathrm{L}(V) \,\middle|\, \exists r_\gamma \in \boldsymbol{F}^\times\ \forall v \in V\colon q(\gamma(v)) = r_\gamma\, \alpha_\gamma(q(v))\right\},$$

$\mathrm{GO}(V,q) \coloneqq \Gamma\mathrm{O}(V,q) \cap \mathrm{GL}(V)$, and $\mathrm{O}(V,q) \coloneqq \left\{\gamma \in \mathrm{GL}(V) \,\middle|\, \forall v \in V\colon q(\gamma(v)) = q(v)\right\}$.

If the form in question is not zero then the *multiplier* $r_\gamma$ of a linear similitude $\gamma$ is determined by $\gamma$, and $\gamma \mapsto r_\gamma$ is a group homomorphism.

The following result dates back to [1, Th. 6, p. 392], see also [5, I § 10, p. 20 with I § 8, p. 15] or [15, Theorem 6.3.1]. For the reader's convenience, we include a proof.

**1.2 Lemma.** *Let $h\colon V \times V \to \boldsymbol{F}$ be a symmetric bilinear form on a vector space $V$ of finite dimension $n$ over a field $\boldsymbol{F}$ with* $\mathrm{char}\,\boldsymbol{F} = 2$. *Let $k \coloneqq \dim_{\boldsymbol{F}^\square} q_h(V)$.*

[1] The treatment in [2] is quite different from that in later editions [4].

**a.** *If $h$ is non-zero then there exists an orthogonal basis for $V$ with respect to $h$ if, and only if, the form $h$ is not alternating, i.e., if there exists $v \in V$ with $q_h(v) \neq 0$.*
**b.** *If $h$ is not alternating then there exists an orthogonal basis $v_1, \dots, v_n$ of $V$ such that $q_h(v_1), \dots, q_h(v_k)$ is an $\boldsymbol{F}^\square$-basis for $q_h(V)$.*

**Proof.** We abbreviate $q \coloneqq q_h$ as in 1.1. If there exists an orthogonal basis then the Gram matrix with respect to that basis is diagonal, and will be zero if the form is alternating. Conversely, assume that there exists $v \in V$ with $q(v) \neq 0$. If $V^\perp \neq \{0\}$, we choose any basis for $V^\perp$, and any vector space complement $W$ to $V^\perp$ with $v \in W$. The restriction of $h$ to $W$ is not degenerate, and not alternating. It suffices to show that there is an orthogonal basis for $W$. Assume that $w_1, \dots, w_k$ are pairwise orthogonal vectors in $W$ with $q(w_i) \neq 0$. Then these vectors are linearly independent, the restriction of $h$ to their span $W_k$ is not degenerate, and $W_k^\perp \cap W$ is a complement to $W_k$ in $W$.

We proceed by induction on $\dim(W_k^\perp \cap W)$: If the restriction of $h$ to $W_k^\perp \cap W$ is either zero or not alternating, we apply the induction hypothesis. It remains to study the case where there exist $x, y \in W_k^\perp \cap W$ with $h(x, y) = 1$ and $q(x) = 0 = q(y)$. Put $w_{k+1} \coloneqq w_k + q(w_k)y$, $w_{k+2} \coloneqq w_k + x + q(w_k)y$, and $\tilde{w}_k \coloneqq w_k + x$. Straightforward computation yields $q(w_{k+1}) = q(w_{k+2}) = q(\tilde{w}_k) = q(w_k) \neq 0$, and that the vectors $\tilde{w}_k, w_{k+1}, w_{k+2} \in \{w_1, \dots, w_{k-1}\}^\perp$ are pairwise orthogonal.

So $w_1, \dots, w_{k-1}, \tilde{w}_k, w_{k+1}, w_{k+2}$ is an orthogonal basis for $W_k + \boldsymbol{F}x + \boldsymbol{F}y$. Applying the induction hypothesis to $(W_k + \boldsymbol{F}x + \boldsymbol{F}y)^\perp \cap W$ finishes the proof of assertion a.

For assertion b, we choose any orthogonal basis and re-order it in such a way that the values of $q$ at the first $k$ basis vectors are linearly independent over $\boldsymbol{F}^\square$. □

**1.3 Lemma.** *Let $h \colon V \times V \to \boldsymbol{F}$ be a symmetric bilinear form over a field $\boldsymbol{F}$ with $\operatorname{char} \boldsymbol{F} = 2$, and let $q \coloneqq q_h \colon V \to \boldsymbol{F} \colon x \mapsto h(x, x)$ be the corresponding quadratic form. If $q$ is anisotropic then the groups $\mathrm{O}(V, q)$ and $\mathrm{O}(V, h)$ are both trivial, and each one of the sets $\boldsymbol{L} \coloneqq \mathrm{GO}(V, q) \cup \{0\}$ and $\boldsymbol{L}_h \coloneqq \mathrm{GO}(V, h) \cup \{0\}$ forms a subfield of the endomorphism ring $\mathrm{End}_{\boldsymbol{F}}(V)$.*

*Putting $r_0 \coloneqq 0$, we extend the multiplier map to a homomorphism $r \colon \boldsymbol{L} \to \boldsymbol{F}$ of fields. This yields isomorphisms from $\boldsymbol{L}$ and from $\boldsymbol{L}_h$, respectively, onto subfields of $\boldsymbol{F}$.*

*The set $q(V)$ is a vector space over the field $r(\boldsymbol{L})$. If $1 \in q(V)$ then $r(\boldsymbol{L}) \subseteq q(V)$.*

**Proof.** Clearly we have $\mathrm{O}(V, h) \leq \mathrm{O}(V, q)$ and $\mathrm{GO}(V, h) \leq \mathrm{GO}(V, q)$. Interpreting $q$ as a semi-linear map from $V$ to the vector space $\boldsymbol{F}$ over $\boldsymbol{F}^\square$ we see that $q$ is injective by our assumption that $q$ is anisotropic. Therefore, the group $\mathrm{O}(V, q)$ is trivial, and so is the subgroup $\mathrm{O}(V, h)$. It follows that the multiplier map $r \colon \mathrm{GO}(V, q) \to \boldsymbol{F}^\times$ is injective. In particular, the group $\mathrm{GO}(V, q)$ is isomorphic to a subgroup of $\boldsymbol{F}^\times$, and thus commutative.

The set $\boldsymbol{L} \coloneqq \mathrm{GO}(V, q) \cup \{0\} \subseteq \mathrm{End}_{\boldsymbol{F}}(V)$ is closed under multiplication, and each member of $\boldsymbol{L} \smallsetminus \{0\}$ has a multiplicative inverse. We extend the multiplier map by putting $r_0 \coloneqq 0$; this extension is still multiplicative.

For $\lambda, \mu \in \boldsymbol{L}$, we note $q\big((\lambda+\mu)(v)\big) = q\big(\lambda(v)+\mu(v)\big) = q\big(\lambda(v)\big) + q\big(\mu(v)\big) = r_\lambda q(v) + r_\mu q(v) = (r_\lambda + r_\mu) q(v)$. So $\boldsymbol{L}$ is closed under addition, and the multiplier map is additive.

Now assume $\lambda, \mu \in \mathrm{GO}(V,h) \le \mathrm{GO}(V,q)$; we need to show $\lambda + \mu \in \boldsymbol{L}_h := \mathrm{GO}(V,h) \cup \{0\}$. The multiplier $(r_\lambda)^2$ of $\lambda^2 \in \mathrm{GO}(V,h)$ is also the multiplier of $r_\lambda \,\mathrm{id} \in \mathrm{GO}(V,h)$, so injectivity of the multiplier map yields $\lambda^2 = r_\lambda \,\mathrm{id}$. Using $\mathrm{char}\, \boldsymbol{F} = 2$ and commutativity of $\mathrm{GO}(V,h)$, we obtain

$$\begin{aligned} h\big(\lambda(v), \mu(w)\big) + h\big(\mu(v), \lambda(w)\big) &= r_\lambda^{-1} h\big(\lambda^2(v), \lambda\mu(w)\big) + r_\mu^{-1} h\big(\mu^2(v), \mu\lambda(w)\big) \\ &= h\big(v, \lambda\mu(w)\big) + h\big(v, \mu\lambda(w)\big) = 0. \end{aligned}$$

This yields $h\big((\lambda+\mu)(v), (\lambda+\mu)(w)\big) = h\big(\lambda(v), \lambda(w)\big) + h\big(\lambda(v), \mu(w)\big) + h\big(\mu(v), \lambda(w)\big) + h\big(\mu(v), \mu(w)\big) = r_\lambda \, h(v,w) + r_\mu \, h(v,w) = (r_\lambda + r_\mu)\, h(v,w)$, as required.

Finally, assume that there exists $v \in V$ with $q(v) = 1$. For $\lambda \in \mathrm{GO}(V,q)$, we obtain $q(\lambda(v)) = r_\lambda q(v) = r_\lambda$, and $q(V) = \{r_\lambda \mid \lambda \in \boldsymbol{L}\}$ follows. □

Clearly, we have $\boldsymbol{F}\mathbf{1} \le \boldsymbol{L}_h \le \boldsymbol{L}$, where $\mathbf{1}$ denotes the identity (matrix). It may well happen that $\boldsymbol{F}\mathbf{1} \ne \boldsymbol{L}_h \ne \boldsymbol{L}$, see 3.4d below.

**1.4 Remark.** Let $\boldsymbol{K} = \boldsymbol{F} + \boldsymbol{F}z$ be the local algebra of degree 2 over a field $\boldsymbol{F}$ with $\mathrm{char}\, \boldsymbol{F} = 2$, with $z^2 = 0$. For a symmetric matrix $S \in \boldsymbol{F}^{3\times3}$, let $g\colon \boldsymbol{K}^3 \times \boldsymbol{K}^3 \to \boldsymbol{K}$ and $\bar{g}\colon \boldsymbol{F}^3 \times \boldsymbol{F}^3 \to \boldsymbol{F}$, respectively, be the bilinear forms with Gram matrix $S$. Every matrix $B \in \boldsymbol{K}^{3\times3}$ can be written as $B = A + zX$ with $A, X \in \boldsymbol{F}^{3\times3}$. Computing $B^\intercal S B = A^\intercal S A + z(A^\intercal S X + X^\intercal S A)$ we see that $B$ belongs to $\mathrm{O}(\boldsymbol{K}^3, g)$ precisely if $A \in \mathrm{O}(\boldsymbol{F}^3, \bar{g})$ and $A^\intercal S X$ is symmetric.

E.g., for $S = \begin{pmatrix} 1&0&0\\0&1&0\\0&0&c \end{pmatrix}$ with $c \in \boldsymbol{F} \smallsetminus \boldsymbol{F}^\square$ we obtain $\mathrm{O}(\boldsymbol{F}^3, \bar{g}) = \left\{ \begin{pmatrix} \mathbf{1}+aJ & 0 \\ 0 & 1 \end{pmatrix} \,\middle|\, a \in \boldsymbol{F} \right\}$, where $J := \begin{pmatrix} 1&1\\1&1 \end{pmatrix}$, and then $\mathrm{O}(\boldsymbol{K}^3, g) = \left\{ \begin{pmatrix} (\mathbf{1}+aJ)(\mathbf{1}+zS) & zc(\mathbf{1}+aJ)v \\ zv^\intercal & 1+zd \end{pmatrix} \,\middle|\, a,d \in \boldsymbol{F}, v \in \boldsymbol{F}^2, S = S^\intercal \in \boldsymbol{F}^{2\times2} \right\}$.

In Section 3, we will repeatedly need the following.

**1.5 Lemma.** *Let $\boldsymbol{K}$ be a commutative local ring with $1+1=0$ in $\boldsymbol{K}$, and let* $\mathrm{SL}_2(\boldsymbol{K}) := \left\{ \begin{pmatrix} a&b\\c&d \end{pmatrix} \,\middle|\, a,b,c,s \in \boldsymbol{K}, ad-bc = 1 \right\}$. *We write* $\boldsymbol{i} := \begin{pmatrix} 0&1\\1&0 \end{pmatrix}$, $L_x := \begin{pmatrix} 1&0\\x&1 \end{pmatrix}$, $U_x := \begin{pmatrix} 1&x\\0&1 \end{pmatrix}$,

$$\hat{L}_x = \begin{pmatrix} 1+x & 0 & x \\ 0 & 1 & 0 \\ x & 0 & 1+x \end{pmatrix}, \quad \hat{U}_x = \begin{pmatrix} 1+x & x & 0 \\ x & 1+x & 0 \\ 0 & 0 & 1 \end{pmatrix}, \quad \text{and} \quad T := \begin{pmatrix} 1&1&1\\1&1&0\\1&0&1 \end{pmatrix}.$$

*The bilinear forms $g$ and $g'$ on $\boldsymbol{K}^3$ are given by $g(x,y) = x_1y_1 + x_2y_3 + x_3y_2$ and $g'(x,y) = \sum_{k=1}^3 x_k y_k$, respectively.*

**a.** *The set $\{L_x \mid x \in \boldsymbol{K}\} \cup \{U_x \mid x \in \boldsymbol{K}\}$ generates* $\mathrm{SL}_2(\boldsymbol{K})$.

**b.** *The set* $\{\hat{L}_x \mid x \in \boldsymbol{K}\} \cup \{\hat{U}_x \mid x \in \boldsymbol{K}\}$ *generates a group* $\hat{\Sigma}$ *isomorphic to* $\mathrm{SL}_2(\boldsymbol{K})$*. Indeed,* $T^{-1}\hat{L}_x T = \left(\begin{smallmatrix} 1 & 0 \\ 0 & L_x \end{smallmatrix}\right)$*,* $T^{-1}\hat{U}_x T = \left(\begin{smallmatrix} 1 & 0 \\ 0 & U_x \end{smallmatrix}\right)$*, and* $T^{-1}\hat{\Sigma}T = \Sigma := \left\{\left(\begin{smallmatrix} 1 & 0 \\ 0 & A \end{smallmatrix}\right) \,\middle|\, A \in \mathrm{SL}_2(\boldsymbol{K})\right\}$*.*

**c.** *Now assume that the maximal ideal* $M$ *of* $\boldsymbol{K}$ *has the property* $MM = \{0\}$*. Then* $\Xi := \left\{\left(\begin{smallmatrix} 1 & u^\intercal \boldsymbol{i} \\ u & \mathbf{1} \end{smallmatrix}\right) \,\middle|\, u \in M^2\right\}$ *and* $\Theta := \left\{\left(\begin{smallmatrix} a & 0 \\ 0 & \mathbf{1} \end{smallmatrix}\right) \,\middle|\, a \in 1 + M\right\}$ *are subgroups, and* $\Xi\Theta$ *is an elementary abelian* 2*-group. The groups* $\Xi$ *and* $\Theta$ *are both normalized by* $\Sigma$*, the semi-direct product*

$$\Xi\Theta\Sigma = \left\{\begin{pmatrix} a & au^\intercal \boldsymbol{i}\, C \\ u & C \end{pmatrix} \,\middle|\, a \in 1 + M,\, u \in M^2,\, C \in \mathrm{SL}_2(\boldsymbol{K})\right\}$$

*equals* $\mathrm{O}(\boldsymbol{K}^3, g)$*, and* $T(\Xi\Theta\Sigma)T^{-1} = \mathrm{O}(V, g')$*.*

*Here* $M^2$ *denotes the set of columns with two entries from* $M$*.*

**Proof.** Assertion a is well known ([7, Th. 4.3.9], see [16, p. 22] for the case where $\boldsymbol{K}$ is a field). Since matrix algebra over a local ring is less popular, we provide a direct argument (actually, in a form that works for every local ring): Let $A = \left(\begin{smallmatrix} a & b \\ c & d \end{smallmatrix}\right) \in \mathrm{SL}_2(\boldsymbol{K})$. Then $a$ and $c$ cannot both be non-invertible. If $a$ is invertible, then

$$\begin{pmatrix} 1 & 0 \\ a^{-1}(c-1) & 1 \end{pmatrix}\begin{pmatrix} 1 & a-1 \\ 0 & 1 \end{pmatrix}\begin{pmatrix} 1 & 0 \\ 1 & 1 \end{pmatrix}\begin{pmatrix} 1 & a^{-1}(1+b-a) \\ 0 & 1 \end{pmatrix} = \begin{pmatrix} a & b \\ c & d \end{pmatrix}.$$

If $c$ is invertible, then

$$\begin{pmatrix} 1 & c^{-1}(a-1) \\ 0 & 1 \end{pmatrix}\begin{pmatrix} 1 & 0 \\ c & 1 \end{pmatrix}\begin{pmatrix} 1 & c^{-1}(d-1) \\ 0 & 1 \end{pmatrix} = \begin{pmatrix} a & b \\ c & d \end{pmatrix}.$$

Assertion b follows by direct computations.

Now assume $MM = \{0\}$. Then $M$ is the kernel of the Frobenius endomorphism, and $\boldsymbol{K}^\square = \{x^2 \mid x \in \boldsymbol{K}\} \cong \boldsymbol{K}/M$ is a field. We interpret the quadratic form $q := q_g$ as a semilinear map from $\boldsymbol{K}^3$ to the vector space $\boldsymbol{K}$ over $\boldsymbol{K}^\square$. The kernel of $q$ is $\{(x_1, x_2, x_3)^\intercal \in \boldsymbol{K}^3 \mid x_1^2 = 0\} = \{(x_1, x_2, x_3)^\intercal \in \boldsymbol{K}^3 \mid x_1 \in M\}$, and contains $\{(0, x_2, x_3)^\intercal \mid x_2, x_3 \in \boldsymbol{K}\}$. That kernel is invariant under each isometry in $\mathrm{O}(V, g)$.

We note $C^\intercal \boldsymbol{i} C = (\det C)\,\boldsymbol{i}$ for each matrix $C \in \boldsymbol{K}^{2\times 2}$ and $v^\intercal \boldsymbol{i}\, v = 0$ for each $v \in \boldsymbol{K}^2$. For any $a \in \boldsymbol{K}$, any $u, w \in \boldsymbol{K}^2$, and $C \in \boldsymbol{K}^{2\times 2}$ we write $w^\intercal = (w_1, w_2)$ and consider the matrix

$$A_{a,u,w,C} := \begin{pmatrix} a & w^\intercal \\ u & C \end{pmatrix}.$$

If the kernel of $q$ is invariant under that matrix, we have $0 = q\big(A_{a,u,w,C}(0,1,0)^\intercal\big) = w_1^2$ and also $0 = q\big(A_{a,u,w,C}(0,0,1)^\intercal\big) = w_2^2$, so $w \in M^2$. Then $ww^\intercal = 0$.

Now $A_{a,u,w,C} \in \mathrm{O}(V,g)$ yields the conditions $a^2 = 1$, $wa = C^\intercal \boldsymbol{i}\, u$, and $C^\intercal \boldsymbol{i}\, C = \boldsymbol{i}$. This leads to $a^{-1} = a \in 1 + M$, $\det C = 1$, and $u \in M^2$.

Conversely, we have $A_{a,u,w,C} \in \mathrm{O}(V,g)$ for each choice of $a \in 1 + M$, $C \in \mathrm{SL}_2(\boldsymbol{K})$, $u \in M^2$ and $w := C^\intercal \boldsymbol{i}\, ua$. So

$$\mathrm{O}(V,g) = \left\{ \begin{pmatrix} a & au^\intercal \boldsymbol{i}\, C \\ u & C \end{pmatrix} \;\middle|\; a \in 1 + M,\, C \in \mathrm{SL}_2(\boldsymbol{K}),\, u \in M^2 \right\}.$$

It is easy to see that both $\Theta := \{A_{a,0,0,\mathbf{1}} \mid a \in 1 + M\}$ and $\Xi$ are elementary abelian 2-groups. A straightforward calculation yields that they centralize each other, and that they are normalized by $\Sigma$. Then $\Theta\Xi$ is an elementary abelian subgroup, and normalized by $\Sigma$, as well. Thus $\mathrm{O}(V,g) = \Xi\Theta\Sigma$ follows.

We note $T^\intercal = T$ and $T^2 = \left(\begin{smallmatrix} 1 & 0 \\ 0 & \boldsymbol{i} \end{smallmatrix}\right)$, so $T^2$ is the Gram matrix for $g$. So $T\,\mathrm{O}(V,g)\,T^{-1}$ equals the group $\mathrm{O}(V,g')$ of isometries with respect to the form with Gram matrix $(T^{-1})^\intercal T^2 T^{-1} = \mathbf{1}$. $\square$

## 2. The Hodge operator in characteristic two

Fundamental definitions and results about Hodge operators have been worked out in some detail in [10]. We repeat the fundamental facts (with simplifications due to the concentration on a special case) and refer to [10] for details and proofs.

**2.1 The Pfaffian form.** Recall that $\dim \bigwedge^n V = 1$ if $\dim V = n$. We fix an isomorphism $b\colon \bigwedge^n V \to \boldsymbol{F}$. For each positive integer $\ell \le n$ the map $b$ induces an isomorphism $\mathrm{Pf}\colon \bigwedge^{n-\ell} V \to \bigwedge^{\ell} V^\vee$ given by

$$\mathrm{Pf}(v_1 \wedge \cdots \wedge v_{n-\ell})(w_1 \wedge \cdots \wedge w_\ell) = b(v_1 \wedge \cdots \wedge v_{n-\ell} \wedge w_1 \wedge \cdots \wedge w_\ell)\,.$$

This is the *Pfaffian form*, see [8, VI 10 Problems 23–28, VIII 12 Problem 42]. As $\operatorname{char} \boldsymbol{F} = 2$, the resulting bilinear map Pf on $\bigwedge V$ is symmetric,

$$\mathrm{Pf}(v_1 \wedge \cdots \wedge v_{n-\ell}, w_1 \wedge \cdots \wedge w_\ell) = \mathrm{Pf}(w_1 \wedge \cdots \wedge w_\ell, v_1 \wedge \cdots \wedge v_{n-\ell})\,.$$

If $n = 2\ell$ is even then $\mathrm{Pf}(v_1 \wedge \cdots \wedge v_\ell, v_1 \wedge \cdots \wedge v_\ell) = 0$ holds for each $v_1 \wedge \cdots \wedge v_\ell$, so Pf is an alternating form on $\bigwedge^\ell V$.

**2.2 Remarks.** For $n = 4$ and $\ell = 2$ we are dealing with the space $\bigwedge^2 \boldsymbol{F}^4$ that carries the Klein quadric. The *quadratic form* Pq defining the Klein quadric is also referred to as a Pfaffian form (cf. [6] and [14] where this form is denoted by $q$), and Pf is the polar form of that quadratic form. Under the present assumption $\operatorname{char} \boldsymbol{F} = 2$, the polar form Pf carries less information than the quadratic form Pq; in fact, for any diagonal quadratic form $k$ (i.e., such that the polar form of $k$ is zero) the sum $\mathrm{Pq} + k$ also has Pf as polar form. Note that Pf is alternating because it is the polar form of a quadratic form in characteristic two.

If one interprets the elements of $\bigwedge^2 \boldsymbol{F}^4$ as alternating matrices then there exists a scalar $s \in \boldsymbol{F}^\times$ such that $\mathrm{Pq}(X)^2 = s \det X$ holds for each $X \in \bigwedge^2 \boldsymbol{F}^4$, cf. [3, § 5 no. 2, Prop. 2, p. 84]; the scalar $s$ reflects the choice of basis underlying that interpretation.

We identify $x_{12} v_1 \wedge v_2 + x_{13} v_1 \wedge v_3 + x_{14} v_1 \wedge v_4 + x_{23} v_2 \wedge v_3 + x_{24} v_2 \wedge v_4 + x_{34} v_3 \wedge v_4$ with

$$X = \begin{pmatrix} 0 & x_{12} & x_{13} & x_{14} \\ -x_{12} & 0 & x_{23} & x_{24} \\ -x_{13} & -x_{23} & 0 & x_{34} \\ -x_{14} & -x_{24} & -x_{34} & 0 \end{pmatrix}$$

and compute $\det X = (x_{12}x_{34} - x_{13}x_{24} + x_{14}x_{23})^2$; up to a scalar (viz., the value $b(v_1 \wedge v_2 \wedge v_3 \wedge v_4)$ above), the value of Pq at $X$ is obtained as $x_{12}x_{34} - x_{13}x_{24} + x_{14}x_{23}$, and the polar form of Pq maps $(X, Y)$ to the corresponding scalar multiple of $x_{12}y_{34} + x_{13}y_{24} + x_{14}y_{23} + x_{23}y_{14} + x_{24}y_{13} + x_{34}y_{12}$.

See [16, 12.14] for an interpretation of Pq in terms of the exterior algebra.

**2.3 The Hodge operator.** We now consider the composite

$$\boldsymbol{J} := \mathrm{Pf}^{-1} \circ \textstyle\bigwedge h \colon \ \bigwedge^\ell V \xrightarrow[\cong]{\bigwedge h} \bigwedge^\ell V^\vee \xrightarrow[\cong]{\mathrm{Pf}^{-1}} \bigwedge^{n-\ell} V \,.$$

This semilinear isomorphism is the *Hodge operator*. It depends, of course, on $h$ and on $b$ but not on the choice of basis.

**2.4 Explicit computation.** Suppose that $v_1, \ldots, v_n$ is an orthogonal basis of $V$. For $\bigwedge^\ell V$ we use the basis vectors $v_{i_1} \wedge \cdots \wedge v_{i_\ell}$ with ascending $i_1 < \cdots < i_\ell \le n$. Then $\bigwedge^\ell h(v_1 \wedge \cdots \wedge v_\ell, -)$ is a linear form on $\bigwedge^\ell V$ which annihilates each one of those basis vectors, except for $v_1 \wedge \cdots \wedge v_\ell$; in fact

$$\textstyle\bigwedge^\ell h(v_1 \wedge \cdots \wedge v_\ell, v_1 \wedge \cdots \wedge v_\ell) = h(v_1, v_1) \cdots h(v_\ell, v_\ell) \,.$$

In other words: $\bigwedge^\ell h$ is again diagonalizable. It then also follows that $\bigwedge^\ell h$ is not degenerate. The linear form $\mathrm{Pf}(v_{\ell+1} \wedge \cdots \wedge v_n)$ annihilates the same collection of basis $\ell$-vectors, and

$$\mathrm{Pf}(v_{\ell+1} \wedge \cdots \wedge v_n, v_1 \wedge \cdots \wedge v_\ell) = b(v_1 \wedge \cdots \wedge v_n) \,.$$

Therefore

$$\boldsymbol{J}(v_1 \wedge \cdots \wedge v_\ell) = v_{\ell+1} \wedge \cdots \wedge v_n \, \frac{h(v_1, v_1) \cdots h(v_\ell, v_\ell)}{b(v_1 \wedge \cdots \wedge v_n)} \,.$$

Note that this last formula is correct only if $v_1, \ldots, v_n$ is an orthogonal basis, and cannot be used if $v_1 \wedge \cdots \wedge v_\ell$ corresponds to a subspace $U$ of $V$ such that the restriction of $h$ to $U$ is degenerate.

**2.5 The square of the Hodge operator.** *Let $H$ be the Gram matrix of $h$ with respect to the orthogonal basis $v_1, \ldots, v_n$. The square of $\boldsymbol{J}$ is a linear automorphism of $\bigwedge^\ell V$, and we find that*

$$\boldsymbol{J}^2 = \delta_\ell \, \mathbf{1} \qquad \textit{where} \qquad \delta_\ell := \tfrac{\det(H)}{b(v_1 \wedge \cdots \wedge v_n)^2} \,.$$

Recall that $\det(H)$ depends on the choice of basis; the invariant would be the square class $\operatorname{disc}(h) \in \boldsymbol{F}^\times / \boldsymbol{F}^\boxtimes$ of $\det(h(v_i, v_j))$. (Here $\boldsymbol{F}^\boxtimes$ denotes the multiplicative group of the field $\boldsymbol{F}^\square := \{s^2 \mid s \in \boldsymbol{F}\}$ of all squares in $\boldsymbol{F}$.) However, the whole expression depends only on $h$ and $b$. Replacing the isomorphism $b\colon \bigwedge^n V \to \boldsymbol{F}$ changes $\boldsymbol{J}$ by a factor and $\boldsymbol{J}^2$ by the square of that factor. In particular, the isomorphism type of the algebra $\boldsymbol{K}_\ell$ introduced in 2.7 below does not depend on the choice of $b$.

**2.6 Lemma.** *For all $x, y \in \bigwedge^\ell V$ we have*

**a.** $\operatorname{Pf}(\boldsymbol{J}(x), y) = \bigwedge^\ell h(x, y)$,

**b.** $\operatorname{Pf}(\boldsymbol{J}(x), \boldsymbol{J}(y)) = \delta_\ell \operatorname{Pf}(y, x)$,

**c.** $\bigwedge^\ell h(\boldsymbol{J}(x), y) = \delta_\ell \operatorname{Pf}(x, y)$,

**d.** $\bigwedge^\ell h(\boldsymbol{J}(x), \boldsymbol{J}(y)) = \delta_\ell \bigwedge^\ell h(x, y)$.

From now on, assume $n = 2\ell$. Then $\boldsymbol{J}$ is an $\boldsymbol{F}$-linear endomorphism of $\bigwedge^\ell V$. We are going to use $\boldsymbol{J}$ to give $\bigwedge^\ell V$ the structure of a right module over an associative algebra of dimension 2 over $\boldsymbol{F}$.

**2.7 The algebra $\boldsymbol{K}_\ell$.** Take $\delta_\ell = \frac{\det(H)}{b(v_1 \wedge \cdots \wedge v_n)^2}$ as in 2.5 and put $\boldsymbol{K}_\ell := \boldsymbol{F}[\boldsymbol{j}_\ell]/(\boldsymbol{j}_\ell^2 - \delta_\ell)$.

**2.8 Definition.** For $v \in \bigwedge^\ell V$ we put $v\,\boldsymbol{j}_\ell := \boldsymbol{J}(v)$. In this way, the space $\bigwedge^\ell V$ becomes an $\mathrm{O}(V, h)$-$\boldsymbol{K}_\ell$-bimodule, i.e., it becomes a right module over $\boldsymbol{K}_\ell$ and $\mathrm{O}(V, h)$ acts $\boldsymbol{K}_\ell$-linearly from the left. Choosing an orthogonal basis $v_1, \ldots, v_n$ for $V$ with a fixed ordering, we obtain a basis $B$ for $\bigwedge^\ell V$ consisting of all $v_{j_1} \wedge \cdots \wedge v_{j_\ell}$ where $(j_1, \ldots, j_\ell)$ is an increasing sequence of length $\ell$ in $\{1, \ldots, n\}$. The sequences with $j_1 = 1$ form a subset $B_1$ of $B$, and $\boldsymbol{J}$ maps each element of $B_1$ to one of $B \smallsetminus B_1$. Moreover, the set $B_1$ forms a basis for the $\boldsymbol{K}_\ell$-module $\bigwedge^\ell V$, showing that the latter is a free module.

**2.9 The bilinear form on the module.** We define $g\colon \bigwedge^\ell V \times \bigwedge^\ell V \to \boldsymbol{K}_\ell$ by

$$g(u, v) := \textstyle\bigwedge^\ell h(u, v) + \bigwedge^\ell h(u, v\boldsymbol{j}_\ell)\, \boldsymbol{j}_\ell^{-1} = \bigwedge^\ell h(u, v) + \boldsymbol{j}_\ell \, (-1)^\ell \operatorname{Pf}(u, v) \,;$$

see 2.6 for the description on the right hand side. This expression is $\boldsymbol{K}_\ell$-bilinear.

Note that $\boldsymbol{K}_\ell$ is not a field, in general: we need the more general concept of bilinear forms over rings.

**2.10 Proposition.** *The form $g$ is diagonalizable.*

For a general proof (and a general formula) see [10, 2.7]. Actually, if $\boldsymbol{K}_\ell$ is a field (of characteristic two, and $\sigma = \mathrm{id}$ by our assumptions) it suffices to note that $g$ is not alternating, see 1.2. For the applications in Section 3 below, we give a special statement explicitly:

**2.11 Example.** Let $v_1, v_2, v_3, v_4$ be an orthogonal basis for $V$, with respect to $h$. Then $w_2 := v_1 \wedge v_2$, $w_3 := v_1 \wedge v_3$, $w_4 := v_1 \wedge v_4$ form an orthogonal basis for the free $\boldsymbol{K}_2$-module $\bigwedge^2 V$, with respect to $g$. Explicitly, we have $g(w_2, w_2) = h(v_1, v_1)h(v_2, v_2)$, $g(w_3, w_3) = h(v_1, v_1)h(v_2, v_2)$, and $g(w_4, w_4) = h(v_1, v_1)h(v_4, v_4)$.

From the definition of $g$ it is clear that $\mathrm{O}(V, h)$ preserves $g$ and that the group $\Gamma\mathrm{O}(V, h)$ of semi-similitudes (see 1.1) acts by semi-similitudes of $g$, see [10, 1.8]. For $\ell = \frac{1}{2}\dim(V)$ we have thus constructed a homomorphism $\eta_\ell \colon \Gamma\mathrm{O}(V, h) \to \Gamma\mathrm{O}(\bigwedge^\ell V, g)$.

**2.12 Lemma.** *The kernel of $\eta_\ell$ is trivial.*

**Proof.** That kernel consists of all scalar multiples $s\mathbf{1}$ of the identity $\mathbf{1}$, where $s^2 = 1$. Since $\operatorname{char} \boldsymbol{F} = 2$, this yields $s = 1$. □

We will call $\boldsymbol{K}_\ell$ *split* whenever it contains divisors of zero. This extends the established terminology for composition algebras. Recall that $\boldsymbol{K}_\ell$ is split precisely if $\delta_\ell$ is a square: $\delta_\ell = s^2$ for some $s \in \boldsymbol{F}^\times$ (and this happens precisely if $h$ has discriminant 1).

In that case, we may assume $s = 1$ without loss of generality. In fact, if we replace our isomorphism $b \colon \bigwedge^n V \to \boldsymbol{F}$ by $sb$ then the Hodge operator $\boldsymbol{J}$ is replaced by $\boldsymbol{J}\, s^{-1} \colon X \mapsto \boldsymbol{J}(X)s^{-1}$, and we have $(\boldsymbol{J}\, s^{-1})^2 = \mathbf{1}$ while the algebra $\boldsymbol{K}_\ell$ remains the same. If $s = 1$ then $z := 1 + \boldsymbol{j}_\ell \in \boldsymbol{K}_\ell$ satisfies $z^2 = 2z = 0$, and is nilpotent. Thus $\boldsymbol{K}_\ell \cong \boldsymbol{F}[X]/(X^2)$ is a local ring if $\boldsymbol{K}_\ell$ is split. Recall that a local ring is a ring in which the set of non-invertible elements is closed under addition (we allow the case where that set consists of 0 alone).

**2.13 Lemma.** *Let $W := \bigwedge^\ell V$, and assume that $\boldsymbol{K}_\ell$ is split.*

- **a.** *The maximal ideal in $\boldsymbol{K}_\ell$ is generated by a nilpotent element $z$. The submodule $Wz$ and the quotient module $W/Wz$ are isomorphic via $\rho_z \colon w + Wz \mapsto wz$.*
- **b.** *The restriction of the form $g$ to the subspace $Wz$ is trivial.*
- **c.** *The $\boldsymbol{K}_\ell$-submodule $Wz$ is invariant under $\eta(\Gamma\mathrm{O}(V, h))$. Thus we obtain a homomorphism $\eta_\ell^o \colon \Gamma\mathrm{O}(V, h) \to \Gamma\mathrm{L}(Wz)$.*
- **d.** *The group induced by $\ker \eta_\ell^o$ on $W$ is an elementary abelian 2-group, acting trivially on $W/Wz$.*

**Proof.** The first three assertions are taken from [10, 3.6]. For the last assertion, we note that elements of $\ker \eta_\ell^o$ act trivially on $W/Wz$ because $\rho_z$ is a module homomorphism. So $\ker \eta_\ell^o$ is isomorphic to a subgroup of $\mathrm{Hom}_{\boldsymbol{F}}(W/Wz, Wz)$. □

The precise structure of $\mathrm{O}(V,h)$ and $\ker \eta_\ell^o$ depends on the defect of the form $h$.

## 3. The four-dimensional cases

We focus on the case where $\ell = 2$ and $n = 2\ell = 4$, and write $\boldsymbol{K} := \boldsymbol{K}_2$. Either $\boldsymbol{K}$ splits and is isomorphic to the local algebra $\boldsymbol{F}[X]/(X^2)$, or we have an inseparable quadratic field extension $\boldsymbol{K}|\boldsymbol{F}$. Recall from 1.1 that $q\colon V \to \boldsymbol{F}\colon x \mapsto h(x,x)$ is a $\varphi$-semilinear map, where $\boldsymbol{F}$ is considered as a vector space over the subfield $\boldsymbol{F}^\square$ of squares, and $\varphi\colon \boldsymbol{F} \to \boldsymbol{F}^\square\colon s \mapsto s^2$ is the Frobenius endomorphism. Note that $\mathrm{O}(V,h)$ is contained in $\mathrm{O}(V,q)$. We distinguish cases according to $\dim_{\boldsymbol{F}^\square} q(V) \in \{1,2,3,4\}$; recall that $\dim_{\boldsymbol{F}} \ker q = 4 - \dim_{\boldsymbol{F}^\square} q(V) \in \{3,2,1,0\}$ is called the defect of $q$. At several places we will use the fact that the orthogonal group of an anisotropic diagonal quadratic form is trivial if the ground field has characteristic 2; cf. [5, §16, p. 35]. Recall from 2.13 that the restriction $g|_{Wz\times Wz}$ is trivial if $z$ is nilpotent (of course, this is of interest only if $\boldsymbol{K}$ splits).

We use the standard basis $e_1, e_2, e_3, e_4$ for $V = \boldsymbol{F}^4$, and write $W := \bigwedge^2 V$.

We will also use the basis $b_1 := (1,1,1,1)^\intercal = e_1+e_2+e_3+e_4$, $b_2 := (1,1,0,0)^\intercal = e_1+e_2$, $b_3 := (1,0,1,0)^\intercal = e_1 + e_3$, $b_4 := (0,0,0,1)^\intercal = e_4$.

**3.1 Proposition.** *If $q$ has defect 3, then $\boldsymbol{K}$ splits and $\mathrm{O}(V,h) \cong \big(\mathrm{SL}_2(\boldsymbol{F}) \ltimes \boldsymbol{F}^2\big) \times \boldsymbol{F}$. The normal subgroup $\Xi \cong \boldsymbol{F}^2 \times \boldsymbol{F}$ is the kernel of the action on $Wz$.*

*Note that $\boldsymbol{K} = \boldsymbol{F} + \boldsymbol{F}z \cong \boldsymbol{F}[X]/(X^2)$ is not a field; the subfield $\boldsymbol{F} < \boldsymbol{K}$ is the image under the endomorphism $\psi\colon \boldsymbol{K} \to \boldsymbol{K}\colon r + sz \mapsto r$. Using matrix descriptions with respect to the $\boldsymbol{K}$-basis $b_1 \wedge b_4$, $b_2 \wedge b_4$, $b_3 \wedge b_4$, we apply $\psi$ to each matrix entry and obtain an endomorphism of the group $\mathrm{GL}_3(\boldsymbol{K})$. Let $\Psi \cong \mathrm{GL}_3(\boldsymbol{F})$ be the image under that endomorphism; in the chosen coordinates, this subgroup $\Psi < \mathrm{GL}_3(\boldsymbol{K})$ consists of all invertible $3 \times 3$ matrices with entries from $\boldsymbol{F}$.*

*We obtain $\mathrm{SL}_2(\boldsymbol{F}) \cong \eta^o(\mathrm{O}(V,h)) = \mathrm{O}(W,g) \cap \Psi < \mathrm{O}(W,g) \cong \mathrm{SL}_2(\boldsymbol{K}) \cong \mathrm{SL}_2(\boldsymbol{F}) \times \mathrm{SL}_2(\boldsymbol{F})$. Every multiplier is a square, the group of similitudes is $\mathrm{GO}(V,h) = \boldsymbol{F}^\times\, \mathrm{O}(V,h)$.*

**Proof.** If $q$ has defect 3 then $\dim_{\boldsymbol{F}^\square} q(V) = 1$ and we may (upon replacing $h$ by a scalar multiple of $h$, cp. 1.2b) assume that $h(x,y) = x^\intercal y$ is the description in coordinates $x$, $y$ with respect to the standard basis. This form $h$ has Witt index 2. In coordinates $v$, $w$ with respect to the basis $b_1, b_2, b_3, b_4$ above, it is given by the map $\tilde h$, with

$$\tilde h(v,w) := v^\intercal \begin{pmatrix} 0 & 0 & 0 & 1 \\ 0 & 0 & 1 & 0 \\ 0 & 1 & 0 & 0 \\ 1 & 0 & 0 & 1 \end{pmatrix} w = v_1w_4 + v_2w_3 + v_3w_2 + v_4w_1 + v_4w_4\,.$$

The orthogonal group $\mathrm{O}(V,\tilde h)$ leaves the quadratic form $\tilde h(v,v) = v_4^2$ invariant. Thus it fixes the linear form with matrix $(0,0,0,1)$, and using suitable block matrices we obtain

$$\mathrm{O}(V,\tilde{h}) \le \left\{ \begin{pmatrix} A & x \\ 0 & a \end{pmatrix} \middle| \; a \in \boldsymbol{F}^{\times}, A \in \mathrm{GL}_3(\boldsymbol{F}), x \in \boldsymbol{F}^3 \right\} .$$

Now one computes easily that

$$\mathrm{O}(V,\tilde{h}) = \left\{ \begin{pmatrix} 1 & t^{\intercal}\boldsymbol{i} & c \\ 0 & B & Bt \\ 0 & 0 & 1 \end{pmatrix} \middle| \; B \in \mathrm{SL}_2(\boldsymbol{F}), t \in \boldsymbol{F}^2, c \in \boldsymbol{F} \right\}$$

where $\boldsymbol{i} := \left(\begin{smallmatrix} 0 & 1 \\ 1 & 0 \end{smallmatrix}\right)$.

Note that the elements with trivial $B$ form an elementary abelian subgroup $\Xi$. In fact, we have an isomorphism

$$\xi \colon \boldsymbol{F}^2 \times \boldsymbol{F} \to \Xi \colon (t_1, t_2, t_3) \mapsto \begin{pmatrix} 1 & t_2 & t_1 & t_3 + t_1 t_2 \\ 0 & 1 & 0 & t_1 \\ 0 & 0 & 1 & t_2 \\ 0 & 0 & 0 & 1 \end{pmatrix} .$$

The map $\xi$ is $\boldsymbol{F}$-linear if we let the scalars act via the rule $(t_1, t_2, t_3) \cdot s := (t_1 s, t_2 s, t_3 s^2)$. In particular, the dimension of $\boldsymbol{F}^2 \times \boldsymbol{F} \cong \Xi$ becomes $2 + \dim_{\boldsymbol{F}^{\square}} \boldsymbol{F}$ which will be greater than 3 whenever the field $\boldsymbol{F}$ is not perfect.[2] The group

$$\Sigma := \left\{ \begin{pmatrix} 1 & 0 & 0 \\ 0 & B & 0 \\ 0 & 0 & 1 \end{pmatrix} \in \mathrm{GL}_4(\boldsymbol{F}) \middle| \; B \in \mathrm{SL}_2(\boldsymbol{F}) \right\} \cong \mathrm{SL}_2(\boldsymbol{F})$$

of block matrices normalizes $\Xi$ and acts in the expected way: it fixes $\xi(\{0\}^2 \times \boldsymbol{F})$ pointwise and induces the usual action on the set

$$\xi\left( \left\{ (t_1, t_2, t_1 t_2) \middle| \; t_1, t_2 \in \boldsymbol{F} \right\} \right) = \left\{ \begin{pmatrix} 1 & t_2 & t_1 & 0 \\ 0 & 1 & 0 & t_1 \\ 0 & 0 & 1 & t_2 \\ 0 & 0 & 0 & 1 \end{pmatrix} \middle| \; t_1, t_2 \in \boldsymbol{F} \right\} .$$

However, that set is not a subgroup; there is no $\Sigma$-invariant subgroup complement to $\xi(\{0\}^2 \times \boldsymbol{F})$.

We fix the isomorphism $b \colon \bigwedge^4 V \to \boldsymbol{F}$ in such a way that $\boldsymbol{J}(e_{\pi(1)} \wedge e_{\pi(2)}) = e_{\pi(3)} \wedge e_{\pi(4)}$ for each permutation $\pi$ of $\{1,2,3,4\}$; recall that the standard basis $e_1, e_2, e_3, e_4$ is an ortho*normal* basis with respect to the form $h$ in the present case. In particular, we now find $\delta = 1$, the algebra $\boldsymbol{K}$ splits, and $z := 1 + \boldsymbol{j}$ is nilpotent.

Using the basis $b_1, \ldots, b_4$ from above, we obtain the $\boldsymbol{K}$-basis $b_1 \wedge b_4$, $b_2 \wedge b_4$, $b_3 \wedge b_4$ for $W$. With respect to that basis, the Gram matrix for $g$ is $\left(\begin{smallmatrix} 1 & 0 & 0 \\ 0 & 0 & 1 \\ 0 & 1 & 0 \end{smallmatrix}\right)$. The elements of $\mathrm{O}(W, g)$ are thus described (with respect to the same basis) by the block matrices $\left(\begin{smallmatrix} a & u^{\intercal}\boldsymbol{i}C \\ u & C \end{smallmatrix}\right)$, with $a \in 1 + Fz$, $u \in (Fz)^2$, and $C \in \mathrm{SL}_2(\boldsymbol{K})$, see 1.5c. We observe that

[2] This phenomenon also plays its role in the study of duality of symplectic quadrangles, cf. [12] and [13].

$$\begin{array}{rclclcl} Y_1 & := & (b_1 \wedge b_4)z & = & (e_1 \wedge e_4 + e_2 \wedge e_4 + e_3 \wedge e_4)z & = & b_1 \wedge b_4 + b_2 \wedge b_3\,, \\ Y_2 & := & (b_2 \wedge b_4)z & = & (e_1 \wedge e_4 + e_2 \wedge e_4)z & = & b_1 \wedge b_2\,, \\ Y_3 & := & (b_3 \wedge b_4)z & = & (e_1 \wedge e_4 + e_3 \wedge e_4)z & = & b_1 \wedge b_3 \end{array}$$

form a basis for $Wz$.

Evaluating $\xi(t_1, t_2, t_3) \in \Xi$ at $Y_1$, $Y_2$, and $Y_3$, we see that $\Xi$ acts trivially on $Wz$, and then also acts trivially on $W/Wz$, see 2.13d. For $B \in \mathrm{SL}_2(\boldsymbol{F})$ we find that $\begin{pmatrix} 1 & 0 & 0 \\ 0 & B & 0 \\ 0 & 0 & 1 \end{pmatrix} \in \Sigma$ maps $Y_1 a_1 + Y_2 a_2 + Y_3 a_3$ to $Y_1 a_1 + Y_2 a_2' + Y_3 a_3'$ with $(a_2', a_3')^\intercal = B(a_2, a_3)^\intercal$. In other words, the image of that element of $\Sigma$ under $\eta^o$ is described by the block matrix $\left(\begin{smallmatrix} 1 & 0 \\ 0 & B \end{smallmatrix}\right)$.

This action of $\mathrm{SL}_2(\boldsymbol{F})$ is an action by isometries of the $\boldsymbol{F}$-bilinear form $k$ on $Wz$ defined by $k(Y_1a_1 + Y_2a_2 + Y_3a_3, Y_1x_1 + Y_2x_2 + Y_3x_3) := a_1x_1 + a_2x_3 + a_3x_2$; see 1.5c (applied with $M = \{0\}$). Note that $\Sigma \cong \mathrm{SL}_2(\boldsymbol{F})$ induces the full group $\mathrm{O}(Wz, k)$. However, the form $k$ is not $g^o$ because $g^o \equiv 0$, see 2.13.

The range $q(V)$ of the quadratic form $q$ is just $\boldsymbol{F}^\square$. So every similitude of $q$ has an element of $\boldsymbol{F}^\boxtimes$ as multiplier, and belongs to $\boldsymbol{F}^\times\, \mathrm{O}(V, q)$. From $\mathrm{GO}(V, h) \le \mathrm{GO}(V, q)$ it then follows that $\mathrm{GO}(V, h) = \boldsymbol{F}^\times\, \mathrm{O}(V, h)$. □

**3.2 Proposition.** *Assume that $q$ has defect 2. In coordinates with respect to the $\boldsymbol{K}$-basis $e_1 \wedge e_2$, $e_1 \wedge e_3$, $e_1 \wedge e_4$, let $\Psi \cong \mathrm{GL}_3(\boldsymbol{F})$ be the subgroup of $\mathrm{GL}_3(\boldsymbol{K})$ consisting of matrices with entries from $\boldsymbol{F}$.*

- **a.** *If $\boldsymbol{K}$ is not split, then $\mathrm{O}(V, h) \cong \mathrm{SL}_2(\boldsymbol{F}) \cong \eta(\mathrm{O}(V, h)) = \mathrm{O}(W, g) \cap \Psi < \mathrm{O}(W, g) \cong \mathrm{SL}_2(\boldsymbol{K})$. The field $\boldsymbol{L}_h$ is an inseparable quadratic extension field of $\boldsymbol{F}\mathbf{1}$.*
- **b.** *If $\boldsymbol{K}$ is split, then $\mathrm{O}(V, h) \cong \boldsymbol{F}^3$ is abelian, and its action on $Wz$ is neither trivial nor faithful. We have $\mathrm{O}(W, g) \cong \boldsymbol{F} \times \big((\boldsymbol{F}^3 \times \boldsymbol{F}^2) \rtimes \boldsymbol{F}\big)$.*
  *Again, the field $\boldsymbol{L}_h$ is an inseparable quadratic extension field of $\boldsymbol{F}\mathbf{1}$.*

*In both cases, we have $\boldsymbol{L} = \boldsymbol{L}_h$, and that field is a quadratic extension of $\boldsymbol{F}\mathbf{1}$.*

**Proof.** If $\dim_{\boldsymbol{F}^\square} q(V) = 2$, we may assume $h(x, y) = x_1c_1y_1 + x_2c_2y_2 + x_3c_3y_3 + x_4c_4y_4$, where $c_3$ and $c_4$ lie in $c_1\boldsymbol{F}^\square + c_2\boldsymbol{F}^\square$, and $c_1, c_2$ are linearly independent over $\boldsymbol{F}^\square$.

If $c_3 \in \boldsymbol{F}^\square c_1$ we may assume $c_3 = c_1$. If $c_3 \notin \boldsymbol{F}^\square c_1$ then there exist $s, t \in \boldsymbol{F}$ with $c_2 = s^2c_1 + t^2c_3$. Then $T := \begin{pmatrix} tc_3 & s \\ sc_1 & t \end{pmatrix}$ is invertible, and $T^\intercal \begin{pmatrix} c_1 & 0 \\ 0 & c_3 \end{pmatrix} T = \begin{pmatrix} c_1^* & 0 \\ 0 & c_2 \end{pmatrix}$ holds with $c_1^* := (tc_3)^2c_1 + (sc_1)^2c_3$. Thus $f_1 := tc_3e_1 + sc_1e_3$, $f_2 := e_2$, $f_3 := se_1 + te_3$, $f_4 := e_4$ is an orthogonal basis, and the Gram matrix for $h$ with respect to that basis has diagonal entries $c_1^*, c_2, c_2, c_4$. Repeating the argument, we obtain that either there exists a diagonal Gram matrix with three identical diagonal entries (if $c_4 \notin \boldsymbol{F}^\square c_1^*$), or there exists a diagonal Gram matrix with two pairs of identical diagonal entries.

Up to similitude, we may thus assume that the Gram matrix (with respect to the standard basis) is one of

$$H_1 := \begin{pmatrix} m & 0 & 0 & 0 \\ 0 & 1 & 0 & 0 \\ 0 & 0 & 1 & 0 \\ 0 & 0 & 0 & 1 \end{pmatrix}, \quad H_2 := \begin{pmatrix} m & 0 & 0 & 0 \\ 0 & m & 0 & 0 \\ 0 & 0 & 1 & 0 \\ 0 & 0 & 0 & 1 \end{pmatrix},$$

respectively, with $m \in \boldsymbol{F} \setminus \boldsymbol{F}^\square$.

**a.** If the form $h$ is described by $H_1$ then its discriminant $m$ is not a square, and $\boldsymbol{K}$ is not split; in fact, we have $\boldsymbol{K} = \boldsymbol{F}(\boldsymbol{j}) \cong \boldsymbol{F}[X]/(X^2 - m)$, with $\boldsymbol{j}^2 = m \notin \boldsymbol{F}^\square$. With respect to the $\boldsymbol{K}$-basis $e_1 \wedge e_2$, $e_1 \wedge e_3$, $e_1 \wedge e_4$, the Gram matrix for $g$ is $\left(\begin{smallmatrix} m & 0 & 0 \\ 0 & m & 0 \\ 0 & 0 & m \end{smallmatrix}\right)$. From 1.5c (applied with $M = \{0\}$) we know that $\mathrm{O}(W, g) \cong \mathrm{SL}_2(\boldsymbol{K})$. As $m$ is a square in $\boldsymbol{K}$, the set $q_g(W)$ coincides with $\boldsymbol{K}^\square$, and every multiplier of $q_g$ is a square. So $\mathrm{GO}(W, q_g) = \boldsymbol{K}^\times \mathrm{O}(W, q_g)$, and $\mathrm{GO}(W, g) = \boldsymbol{K}^\times \mathrm{O}(W, g)$.

In order to understand the group $\mathrm{O}(V, h)$, we first study the quadratic form given by $q(x) = h(x, x) = mx_1^2 + x_2^2 + x_3^2 + x_4^2 = mx_1^2 + (x_2 + x_3 + x_4)^2$. As $m$ is not a square in $\boldsymbol{F}$, the kernel of $q$ is the hyperplane $\left\{ (0, x_2, x_3, x_4)^\intercal \,\middle|\, x_2 + x_3 + x_4 = 0 \right\}$. This hyperplane is invariant under $\mathrm{O}(V, h)$; it will be convenient to use the basis $b_1 := e_1$, $b_2 := e_2 + e_3 + e_4$, $b_3 := e_2 + e_3$, $b_4 := e_2 + e_4$. With respect to that basis, the Gram matrix of $h$ is the block matrix $\left(\begin{smallmatrix} N & 0 \\ 0 & \boldsymbol{i} \end{smallmatrix}\right)$, where $N := \left(\begin{smallmatrix} m & 0 \\ 0 & 1 \end{smallmatrix}\right)$ is a Gram matrix for the norm form of $\boldsymbol{K}|\boldsymbol{F}$, and $\boldsymbol{i} := \left(\begin{smallmatrix} 0 & 1 \\ 1 & 0 \end{smallmatrix}\right)$. In coordinates with respect to that basis, the isometry group of $h$ consists of the block matrices of the form $\left(\begin{smallmatrix} E & 0 \\ 0 & A \end{smallmatrix}\right)$ with the $2 \times 2$ identity matrix $E$, and $A \in \mathrm{SL}_2(\boldsymbol{F})$. In particular, we find $\mathrm{O}(V, h) \cong \mathrm{SL}_2(\boldsymbol{F})$.

As in 1.5a, we generate the group $\mathrm{SL}_2(\boldsymbol{F})$ by the matrices $L_x$ and $U_x$, with $x \in \boldsymbol{F}$. Transforming $\left(\begin{smallmatrix} E & 0 \\ 0 & L_x \end{smallmatrix}\right)$ and $\left(\begin{smallmatrix} E & 0 \\ 0 & U_x \end{smallmatrix}\right)$ back into the description with respect to standard coordinates, we obtain that $\mathrm{O}(V, h)$ is generated by the matrices $\tilde{A}_x := \tilde{T}^{-1} \left(\begin{smallmatrix} E & 0 \\ 0 & A_x \end{smallmatrix}\right) \tilde{T}$ and $\tilde{B}_x := \tilde{T}^{-1} \left(\begin{smallmatrix} E & 0 \\ 0 & B_x \end{smallmatrix}\right) \tilde{T}$, where $\tilde{T} = \left(\begin{smallmatrix} 1 & 0 \\ 0 & T \end{smallmatrix}\right)$, $\tilde{A}_x = \left(\begin{smallmatrix} 1 & 0 \\ 0 & \hat{L}_x \end{smallmatrix}\right)$, and $\tilde{B}_x = \left(\begin{smallmatrix} 1 & 0 \\ 0 & \hat{U}_x \end{smallmatrix}\right)$, with $x \in \boldsymbol{F}$.

With respect to the $\boldsymbol{K}$-basis $e_1 \wedge e_2$, $e_1 \wedge e_3$, $e_1 \wedge e_4$ for $W$, we then find that the action of these elements on $W$ is described by the matrices $\hat{L}_x$ and $\hat{U}_x$, respectively, with $x \in \boldsymbol{F}$. The same matrices, but with $x \in \boldsymbol{K}$ instead of $x \in \boldsymbol{F}$, generate $\mathrm{O}(W, g) \cong \mathrm{SL}_2(\boldsymbol{K})$; see 1.5c. Therefore, the image of $\mathrm{O}(V, h)$ under $\eta$ equals $\mathrm{O}(W, g) \cap \mathrm{GL}_3(\boldsymbol{F})$.

For each similitude $\varphi \in \mathrm{GO}(V, h)$, the multiplier $r_\varphi$ lies in the range $q(V) = \boldsymbol{F}^\square + \boldsymbol{F}^\square m$ because that range contains 1. We note that $q(V)$ forms the quadratic extension field $\boldsymbol{F}^\square(m) \cong \boldsymbol{F}(\sqrt{m}) \cong \boldsymbol{K}$ of $\boldsymbol{F}^\square$. Every element $a^2 + b^2 m \in \boldsymbol{F}^\square(m) \setminus \{0\}$ is the multiplier of some similitude of the form $h$; in coordinates with respect to the basis $b_1, b_2, b_3, b_4$, the block matrix $\left(\begin{smallmatrix} A & 0 \\ 0 & A \end{smallmatrix}\right)$ with $A := \left(\begin{smallmatrix} a & b \\ bm & a \end{smallmatrix}\right) \in \mathrm{GL}_2(\boldsymbol{F})$ describes a similitude with multiplier $a^2 + b^2 m$. So $\boldsymbol{L} = \boldsymbol{L}_h = \boldsymbol{F}\mathbf{1} + \boldsymbol{F} \left(\begin{smallmatrix} 0 & 1 \\ m & 0 \end{smallmatrix}\right)$.

**b.** Now consider the case where $h$ is described by $H_2$. Then the discriminant is a square, and $\boldsymbol{K}$ is split. We normalize such that $\boldsymbol{j}^2 = 1$, and $\boldsymbol{K} = \boldsymbol{F} + \boldsymbol{F}z \cong \boldsymbol{F}[X]/(X^2)$ holds for $z := 1 + \boldsymbol{j}$. The kernel of the quadratic form given by $q(x) = h(x, x)$ is spanned by $d_1 := e_1 + e_2$ and $d_2 := e_3 + e_4$. We use $d_3 := e_1$ and $d_4 := e_4$ to extend this to a basis for $V$. With respect to that basis, the Gram matrix for $h$ is the block matrix $\left(\begin{smallmatrix} 0 & N \\ N & N \end{smallmatrix}\right)$. Again, $N = \left(\begin{smallmatrix} m & 0 \\ 0 & 1 \end{smallmatrix}\right)$ but this is no longer a Gram matrix for the norm of $\boldsymbol{K}|\boldsymbol{F}$. However, from $2 = \dim_{\boldsymbol{F}^\square} q(V) = \dim(\boldsymbol{F}^\square + \boldsymbol{F}^\square m)$ we infer that $1, m$ are linearly independent over $\boldsymbol{F}^\square$, and $N$ is the Gram matrix of an anisotropic quadratic form on $\boldsymbol{F}^2$. Thus the isometry group of that form is trivial.

In the chosen coordinates, a block matrix $\left(\begin{smallmatrix} A & B \\ C & D \end{smallmatrix}\right)$ describing an isometry of $h$ will leave the kernel of $q$ invariant, so $C = 0$. Invariance of $q$ then gives $\left(\begin{smallmatrix} 0 & 0 \\ 0 & N \end{smallmatrix}\right) = \left(\begin{smallmatrix} A^\intercal & 0 \\ B^\intercal & D^\intercal \end{smallmatrix}\right)\left(\begin{smallmatrix} 0 & 0 \\ 0 & N \end{smallmatrix}\right)\left(\begin{smallmatrix} A & B \\ 0 & D \end{smallmatrix}\right) = \left(\begin{smallmatrix} 0 & 0 \\ 0 & D^\intercal N D \end{smallmatrix}\right)$, and $D$ is an element of the trivial isometry group of the quadratic form described by $N$. So $D = \mathbf{1}$. Now invariance of $h$ implies $\left(\begin{smallmatrix} 0 & N \\ N & N \end{smallmatrix}\right) = \left(\begin{smallmatrix} A^\intercal & 0 \\ B^\intercal & \mathbf{1} \end{smallmatrix}\right)\left(\begin{smallmatrix} 0 & N \\ N & N \end{smallmatrix}\right)\left(\begin{smallmatrix} A & B \\ 0 & \mathbf{1} \end{smallmatrix}\right) = \left(\begin{smallmatrix} 0 & A^\intercal N \\ NA & NB+B^\intercal N+N \end{smallmatrix}\right)$, and $A = \mathbf{1}$ follows.

In the chosen coordinates, the isometry group of $h$ therefore consists of the block matrices of the form $\left(\begin{smallmatrix} \mathbf{1} & B \\ 0 & \mathbf{1} \end{smallmatrix}\right)$ with $B \in \boldsymbol{F}^{2\times 2}$ such that $NB$ is a symmetric matrix. Transforming this description into coordinates with respect to the standard basis $e_1, e_2, e_3, e_4$, we obtain $\mathrm{O}(V, h) = \left\{ \left(\begin{smallmatrix} \mathbf{1}+aJ & bJ \\ mbJ & \mathbf{1}+cJ \end{smallmatrix}\right) \,\middle|\, a, b, c \in F \right\}$, where $J = \left(\begin{smallmatrix} 1 & 1 \\ 1 & 1 \end{smallmatrix}\right)$. This shows that $\mathrm{O}(V, h) \cong \boldsymbol{F}^3$ is an elementary abelian 2-group.

Now $w_1 := e_1 \wedge e_4$, $w_2 := e_1 \wedge e_3$, $w_3 := e_1 \wedge e_2$ is a $\boldsymbol{K}$-basis for $W$, with $\boldsymbol{j}w_1 = e_2 \wedge e_3$, $\boldsymbol{j}w_2 = e_2 \wedge e_4$, $\boldsymbol{j}w_3 = me_3 \wedge e_4$. With respect to the $K$-basis $w_1, w_2, w_3$, the action of $\left(\begin{smallmatrix} \mathbf{1}+aJ & bJ \\ mbJ & \mathbf{1}+cJ \end{smallmatrix}\right) \in \mathrm{O}(V, h)$ is described by the matrix

$$\begin{pmatrix} 1+a+c+z(ac+mb^2) & a+c+z(a+ac+mb^2) & zmb \\ a+c+z(a+ac+mb^2) & 1+a+c+z(ac+mb^2) & zmb \\ zb & zb & 1 \end{pmatrix},$$

and the action on $Wz$, with respect to the $\boldsymbol{F}$-basis

$$\begin{aligned} w_1 z &:= (e_1 \wedge e_4)z = e_1 \wedge e_4 + e_2 \wedge e_3, \\ w_2 z &:= (e_1 \wedge e_3)z = e_1 \wedge e_3 + e_2 \wedge e_4, \\ w_3 z &:= (e_1 \wedge e_2)z = e_1 \wedge e_2 + e_3 \wedge e_4, \end{aligned}$$

is described by the matrix

$$\begin{pmatrix} 1+a+c & a+c & 0 \\ a+c & 1+a+c & 0 \\ 0 & 0 & 1 \end{pmatrix}.$$

This shows that the action of $\mathrm{O}(V, h)$ on $Wz$ is not faithful, the kernel of the restriction $\eta^o$ is

$$\ker \eta^o = \left\{ \begin{pmatrix} \mathbf{1}+aJ & bJ \\ bJ & \mathbf{1}+aJ \end{pmatrix} \,\middle|\, a \in \boldsymbol{F} \right\}.$$

In coordinates with respect to the basis $w_1, w_3, w_1 + w_3$, the form $g$ is described by the matrix $\left(\begin{smallmatrix} m & 0 & m \\ 0 & 1 & 0 \\ m & 0 & 0 \end{smallmatrix}\right)$. The matrix describing an isometry $\lambda$ of $g$ with respect to that basis has all entries from $F$ if, and only if, that matrix is in $\left\{ \left( \begin{smallmatrix} 1&0&0 \\ 0&1&0 \\ a&0&1 \end{smallmatrix} \right) \middle| \; a \in \boldsymbol{F} \right\}$. With respect to the basis $w_1, w_2, w_3$, such an isometry $\lambda$ is thus described by a matrix of the form $\left( \begin{smallmatrix} \mathbf{1}+aJ & 0 \\ 0 & 1 \end{smallmatrix} \right)$. Using 1.4 we now infer that the elements of $\mathrm{O}(W,g)$ are exactly those $K$-linear maps that are described (with respect to $w_1, w_2, w_3$) by matrices of the form

$$\psi(a,u,S,t) := \begin{pmatrix} (\mathbf{1}+aJ)(\mathbf{1}+zS) & zm(\mathbf{1}+aJ)u \\ zu^{\intercal} & 1+zt \end{pmatrix},$$

where $a, t \in \boldsymbol{F}$, $u \in \boldsymbol{F}^2$, and $S$ is a symmetric matrix with entries from $\boldsymbol{F}$. So

$$\mathrm{O}(W,g) = \left\{ \psi(a,u,S,t) \,\middle|\, a,t \in F,\ u \in F^2,\ S = S^{\intercal} \in F^{2\times 2} \right\};$$

the multiplication is given by

$$\psi(a,u,S,s)\,\psi(b,v,T,t) = \psi\left(a+b, u+v+bJu, S+T+\mathrm{tr}(S)b\boldsymbol{i}, s+t\right).$$

The subgroups $\Theta := \{\psi(a,0,0,0) \mid a \in \boldsymbol{F}\}$ and $\mathrm{Z} := \{\psi(0,0,0,s) \mid s \in \boldsymbol{F}\}$ are both in a natural way isomorphic to the additive group of $\boldsymbol{F}$. Note that $\mathrm{Z}$ is contained in the center

$$\left\{ \psi(0,u,S,s) \,\middle|\, s \in \boldsymbol{F}, u \in JF^2, S = S^{\intercal} \in \boldsymbol{F}^{3\times 3}, \mathrm{tr}(S) = 0 \right\}$$

of $\mathrm{O}(W,g)$. The subgroups

$$\Xi_2 := \{\psi(0,u,0,0) \mid u \in F^2\} \cong \boldsymbol{F}^2 \quad \text{and} \quad \Xi_3 := \{\psi(0,0,S,0) \mid S = S^{\intercal} \in \boldsymbol{F}^{2\times 2}\} \cong \boldsymbol{F}^3$$

are both normal in $\mathrm{O}(W,g)$, and $\mathrm{O}(W,g)$ is the direct product of $\mathrm{Z}$ with the semidirect product of $\Xi_2 \times \Xi_3$ with the subgroup $\Theta$; conjugation by $\psi(a,0,0,0)$ maps $\psi(0,u,S,0) \in \Xi_2\,\Xi_3$ to $\psi\big(0,\ u+aJu,\ S+\mathrm{tr}(S)(a\boldsymbol{i}+a^2J), 0\big)$.

Note that $\eta$ maps $\left( \begin{smallmatrix} \mathbf{1}+aJ & bJ \\ mbJ & \mathbf{1}+cJ \end{smallmatrix} \right)$ to $\psi\left(a+c, \left(\begin{smallmatrix} b \\ b \end{smallmatrix}\right), \left( \begin{smallmatrix} a^2+b^2m & a+a^2+b^2m \\ a+a^2+b^2m & a^2+b^2m \end{smallmatrix} \right), 0\right)$.

As in case a, the set $q(V)$ coincides with the subfield $\boldsymbol{F}^{\square}(m)$ of $\boldsymbol{F}$, and $\boldsymbol{F}^{\square}(m) \setminus \{0\}$ is the set of all multipliers of similitudes of $h$; in coordinates with respect to the basis $d_1, d_2, d_3, d_4$, the block matrix $\left( \begin{smallmatrix} A & 0 \\ 0 & A \end{smallmatrix} \right)$ with $A := \left( \begin{smallmatrix} a & b \\ bm & a \end{smallmatrix} \right) \in \mathrm{GL}_2(\boldsymbol{F})$ describes a similitude with multiplier $a^2 + b^2m$. So $\boldsymbol{L}_h = \boldsymbol{L} = \boldsymbol{F}\mathbf{1} + \boldsymbol{F}\left( \begin{smallmatrix} 0 & 1 \\ m & 0 \end{smallmatrix} \right)$ is an extension of degree 2 over $\boldsymbol{F}\mathbf{1}$. □

**3.3 Proposition.** *If $q$ has defect 1, then $\mathrm{O}(V,h) \cong \boldsymbol{F}$ is abelian, and $\boldsymbol{K}$ is not split. We have $\mathrm{O}(V,h) \cong (\boldsymbol{F},+)$, and $\eta(\mathrm{O}(V,h)) < \mathrm{O}(W,g) \cong (\boldsymbol{K},+)$. There is a basis $v_1, v_2, v_3, v_4$ for $V$ such that $h$ has the Gram matrix*

$$H := \begin{pmatrix} 0 & 1 & 0 & 0 \\ 1 & 1 & 0 & 0 \\ 0 & 0 & c_3 & 0 \\ 0 & 0 & 0 & c_4 \end{pmatrix}.$$

*With respect to the basis* $v_1 \wedge v_3$, $v_2 \wedge v_3$, $v_1 \wedge v_2$ *for* $W$, *the elements of* $\eta(\mathrm{O}(V,h))$ *are those elements of* $\mathrm{O}(W,g)$ *that are described by matrices with entries from* $\boldsymbol{F}$.

*For each one of the forms* $h$, $q_h$, *and* $g$, *every multiplier is a square, the respective groups of similitudes are* $\mathrm{GO}(V,h) = \boldsymbol{F}^\times \mathrm{O}(V,h)$ *and* $\mathrm{GO}(W,g) = \boldsymbol{K}^\times \mathrm{O}(W,g)$. *However, the multipliers of the form* $q_g$ *are the non-zero elements of the quadratic extension* $\boldsymbol{K}^\square(c_3)$ *over* $\boldsymbol{K}^\square$.

**Proof.** If $\dim_{\boldsymbol{F}^\square} q(V) = 3$ then there is $u_1 \in V \smallsetminus \{0\}$ with $h(u_1,u_1) = 0$. As $h$ is not degenerate, there exists $u_2 \in V$ with $h(u_1,u_2) = 1$. We write $s := h(u_2,u_2)$. If $s = 0$ then $\dim\ker q > 1$, contradicting our assumption that the range of $q$ has dimension 3. We multiply the form by $s^{-1}$ and then replace $u_1$ by $su_1$; we may thus assume $s = 1$.

The restriction of $h$ to $\boldsymbol{F}u_1 + \boldsymbol{F}u_2$ is not degenerate, so $\{u_1,u_2\}^\perp$ forms a vector space complement for that subspace in $V$. If the restriction of $h$ to that complement were isotropic then $\dim\ker q$ would be greater than 1.

So the restriction of $h$ to $\{u_1,u_2\}^\perp$ is anisotropic, and diagonalizable by 1.2. Choosing an orthonormal basis $u_3, u_4$ for $\{u_1,u_2\}^\perp$, we obtain that the Gram matrix for $h$ with respect to the basis $u_1,u_2,u_3,u_4$ is $H$; the discriminant of $h$ is represented by $\delta := c_3c_4$. Now $\delta \notin \boldsymbol{F}^\square$; otherwise, the vector $u_3c_4 + u_4\sqrt{c_3c_4}$ would be isotropic, and $\dim\ker q$ would be greater than 1. So $\boldsymbol{K} \cong \boldsymbol{F}(\sqrt{\delta})$ is not split. We note that $c_3 \notin \boldsymbol{K}^\square$; in fact $c_3 = \left(r + t\sqrt{\delta}\right)^2$ with $r,t \in \boldsymbol{F}$ would imply $0 = 1 + c_3(r/c_3)^2 + c_4t^2$, contradicting the fact that $\dim\ker q = 1$.

The group $\mathrm{O}(V,h)$ leaves invariant the quadratic form $q$. Therefore the subspace $\ker q$ is invariant under $\mathrm{O}(V,h)$, and so is the orthogonal space $\ker q^\perp$, which is spanned by $\{u_1,u_3,u_4\}$. Using these facts facilitates to see that the elements of $\mathrm{O}(V,h)$ are described (with respect to the basis $u_1,u_2,u_3,u_4$) by matrices of the form $\left(\begin{smallmatrix} U_x & 0 \\ 0 & \mathbf{1} \end{smallmatrix}\right)$ with $x \in \boldsymbol{F}$, $E = \left(\begin{smallmatrix} 1 & 0 \\ 0 & 1 \end{smallmatrix}\right)$, and $U_x$ as in 1.5a. So $\mathrm{O}(V,h) \cong (\boldsymbol{F},+)$.

The vectors $v_1 := u_1 + u_2$, $v_2 := u_2$, $v_3 := u_3$, $v_4 := u_4$ form an orthogonal basis. In coordinates with respect to that basis, the group $\mathrm{O}(V,h)$ consists of the block matrices $\tilde{U}_x := \left(\begin{smallmatrix} \hat{U}_x & 0 \\ 0 & 1 \end{smallmatrix}\right)$, with $x \in \boldsymbol{F}$. With respect to the $\boldsymbol{K}$-basis $w_1 := v_1 \wedge v_3$, $w_2 := v_2 \wedge v_3$, $w_3 := v_1 \wedge v_2$, the Gram matrix for the bilinear form $g$ is $\left(\begin{smallmatrix} c_3 & 0 & 0 \\ 0 & c_3 & 0 \\ 0 & 0 & 1 \end{smallmatrix}\right)$. The form $g$ is isotropic (as prophesied by [10, 2.9]); in fact, the vector $w_1 + w_2$ is isotropic with respect to $g$.

In order to determine $\mathrm{GO}(W,g)$, we consider the quadratic form $q_g$ obtained by evaluating $g$ on the diagonal. That quadratic form has one-dimensional kernel because $c_3 \notin \boldsymbol{K}^\square$. So $\boldsymbol{K}(w_1 + w_2)$ is that kernel, and both $\boldsymbol{K}(w_1 + w_2)$ and $(w_1 + w_2)^\perp = \boldsymbol{K}(w_1 + w_2) + \boldsymbol{K}w_3$ are invariant under the group $\mathrm{GO}(W,g)$. It turns out that the elements of $\mathrm{GO}(W,g)$ are described (with respect to the basis $w_1,w_2,w_3$) by matrices of the form $a\hat{U}_x$, with $a \in \boldsymbol{K}^\times$ and $\hat{U}_x$ as in 1.5, with $x \in \boldsymbol{K}$. This yields $\mathrm{O}(W,g) \cong (\boldsymbol{K},+)$, and the group $\mathrm{GO}(W,g)$ is the direct product of that group by the multiplicative group $\boldsymbol{K}^\times$. In coordinates with respect to the basis $w_1 + w_2, w_2, w_3$, the members of $\mathrm{O}(W,q_g)$ are described by matrices of the form $\left(\begin{smallmatrix}1 & b\\ 0 & \mathbf{1}\end{smallmatrix}\right)$ with $b \in K^2$, and $\begin{pmatrix}1&0&0\\0&0&c_3\\0&1&0\end{pmatrix}$ describes a similitude of $q_g$ with multiplier $c_3$, so $\boldsymbol{K}^\square(c_3)$ is the field of multipliers for $q_g$.

For each similitude of $h$, the multiplier lies in $q(V) = \boldsymbol{F}^\square + \boldsymbol{F}^\square c_3 + \boldsymbol{F}^\square c_4$. If $m \in q(V) \smallsetminus \boldsymbol{F}^\square$ were a multiplier of some similitude then $q(V)$ would be a vector space over the field $\boldsymbol{F}^\square(m)$. Since this is impossible, we obtain $\mathrm{GO}(V,h) = \boldsymbol{F}^\times \mathrm{O}(V,h)$, and $\mathrm{GO}(V,q) = \boldsymbol{F}^\times \mathrm{O}(V,q)$. □

**3.4 Proposition.** *If $q$ has defect 0 then both $\mathrm{O}(V,q)$ and $\mathrm{O}(V,h)$ are trivial, and the following assertions hold.*

- **a.** *If $q(V)$ is a subfield of $\boldsymbol{F}$ then the subfield $\boldsymbol{L} := \mathrm{GO}(V,q) \cup \{0\}$ of $\mathrm{End}_{\boldsymbol{F}}(V)$ is a totally inseparable extension of degree 4 and exponent 1 over $\boldsymbol{F}\mathbf{1}$, and $r(\boldsymbol{L}) = q(V)$. In all other cases, we have $\boldsymbol{L} = \boldsymbol{F}\mathbf{1}$.*
- **b.** *If there exists $\delta \in \mathrm{GO}(V,h) \smallsetminus \boldsymbol{F}\mathbf{1}$ then the discriminant of $h$ is represented by an element of the field $\boldsymbol{F}^\square(r_\delta)$.*
- **c.** *If the discriminant of $h$ is trivial then the field $\boldsymbol{L}_h := \mathrm{GO}(V,h) \cup \{0\}$ coincides with $\boldsymbol{L}$, and $r(\boldsymbol{L}_h) = r(\boldsymbol{L}) = q(V)$.*
- **d.** *If the discriminant of $h$ is not trivial but $q(V)$ is a subfield of $\boldsymbol{F}$ then the field $\boldsymbol{L}_h$ is an inseparable extension of degree 2 over $\boldsymbol{F}\mathbf{1}$, and $\boldsymbol{F}^\square < r(\boldsymbol{L}_h) < r(\boldsymbol{L}) = q(V)$.*
- **e.** *If $1 \in q(V)$ but $q(V)$ is not a subfield of $\boldsymbol{F}$ then $\boldsymbol{L}_h = \boldsymbol{L} = \boldsymbol{F}\mathbf{1}$, and the discriminant of $h$ is not trivial.*

**Proof.** Let $v_1, v_2, v_3, v_4$ be an orthogonal basis for $V$ with respect to $h$. Passing to a scalar multiple of the form, we may assume $q(v_1) = 1$, then both $h$ and $q$ are described by the diagonal matrix with diagonal entries $1$, $a := q(v_2)$, $b := q(v_3)$, and $c := q(v_4)$. The discriminant of $h$ is the coset in $\boldsymbol{F}^\times/\boldsymbol{F}^\boxtimes$ that is represented by $\det H = abc$, and $1, a, b, c$ are linearly independent over $\boldsymbol{F}^\square$ because $q$ is anisotropic.

Recall from 1.3 that $\boldsymbol{L}$ and $\boldsymbol{L}_h$ are subfields of the endomorphism algebra $\mathrm{End}_{\boldsymbol{F}}(V)$, and that the multiplier map $r\colon \boldsymbol{L} \to \boldsymbol{F}$ is a field homomorphism. For $\lambda \in \mathrm{GO}(V,q)$, the set $\boldsymbol{E}_\lambda := \boldsymbol{F}\mathbf{1} + \boldsymbol{F}\lambda$ is a subfield of $\boldsymbol{L}$, and $q(V)$ is a vector space over the field $r(\boldsymbol{E}_\lambda) = \{r_\gamma \mid \gamma \in \boldsymbol{F}\mathbf{1} + \boldsymbol{F}\lambda\} = \boldsymbol{F}^\square(r_\lambda)$.

If $\boldsymbol{L} \neq \boldsymbol{F}\mathbf{1}$, we choose $\lambda \in \mathrm{GO}(V,q)$ such that the multiplier $r_\lambda$ is not a square. Then $q(V)$ has dimension 2 over $r(\boldsymbol{E}_\lambda)$. The restriction of $h$ to $\boldsymbol{E}_\lambda v_1 = \boldsymbol{F}v_1 + \boldsymbol{F}\lambda(v_1)$ is not degenerate because (the restriction of) $q$ is anisotropic.

Pick any $w \in (\boldsymbol{E}_\lambda v_1)^\perp \smallsetminus \{0\}$. Then $d := q(w)$ is not in $r(\boldsymbol{E}_\lambda) = q(\boldsymbol{E}_\lambda(v_1))$ because $q$ is injective. Now $q(V) = q(\boldsymbol{E}_\lambda) + q(\boldsymbol{E}_\lambda)d = \boldsymbol{F}^\square + \boldsymbol{F}^\square r_\lambda + \boldsymbol{F}^\square d + \boldsymbol{F}^\square r_\lambda d$ coincides with the field $\boldsymbol{F}^\square(r_\lambda, d)$ generated by $r_\lambda$ and $d$ over $\boldsymbol{F}^\square$. This is a totally inseparable extension of degree 4 and exponent 1.

Conversely, assume that $q(V)$ is a subfield of $\boldsymbol{F}$. Then there are $s_1, s_a, s_b, s_c \in \boldsymbol{F}$ with $ab = s_1^2 + s_a^2 a + s_b^2 b + s_c^2 c$, and $s_c \neq 0$ because $1, a, b, ab$ are linearly independent over $\boldsymbol{F}^\square$. Replacing the fourth basis vector $v_4$ by $s_1 v_1 + s_a v_2 + s_b v_3 + s_c v_4$ now yields a basis for $V$, and with respect to that basis, the form $q$ is described by the diagonal matrix with diagonal entries $1, a, b, ab$. In those coordinates, the matrices

$$A := \begin{pmatrix} 0 & a & 0 & 0 \\ 1 & 0 & 0 & 0 \\ 0 & 0 & 0 & a \\ 0 & 0 & 1 & 0 \end{pmatrix} \text{ and } B := \begin{pmatrix} 0 & 0 & b & 0 \\ 0 & 0 & 0 & b \\ 1 & 0 & 0 & 0 \\ 0 & 1 & 0 & 0 \end{pmatrix}$$

represent similitudes of $q$, and taking multipliers maps $\boldsymbol{L} = \boldsymbol{F}\boldsymbol{1} + \boldsymbol{F}A + \boldsymbol{F}B + \boldsymbol{F}AB$ onto $\boldsymbol{F}^{\square}(a, b) = q(V)$. This completes the proof of assertion a.

Now assume $\lambda \in \mathrm{GO}(V, h)$, and that $\lambda$ is not a scalar multiple of the identity. Then $q(V)$ is a subfield of $\boldsymbol{F}$, and $r_\lambda \in q(V) \smallsetminus \boldsymbol{F}^{\square}$. We put $w_1 := v_1$, $w_2 := \lambda(v_1)$, and pick $w_3 \in \{w_1, w_2\}^{\perp} \smallsetminus \{0\}$. We observe that the restriction of the form $h$ to $\boldsymbol{F}w_1 + \boldsymbol{F}w_2$ is not degenerate. As $\lambda$ is a similitude of $h$, the subspace $\boldsymbol{F}w_3 + \boldsymbol{F}\lambda(w_3)$ is the orthogonal complement for $\boldsymbol{F}w_1 + \boldsymbol{F}w_2$.

With respect to the basis $w_1, w_2, w_3, w_4 := \lambda(v_3)$, the form $h$ is described by the Gram matrix

$$H' := \begin{pmatrix} 1 & s & 0 & 0 \\ s & r_\lambda & 0 & 0 \\ 0 & 0 & q(w_3) & t \\ 0 & 0 & t & r_\lambda q(w_3) \end{pmatrix},$$

where $s := h(w_1, w_2)$, and $t = h(w_3, w_4)$. The discriminant of $h$ is represented by $\det H' = (r_\lambda + s^2)(q(w)^2 r_\lambda + t^2) = q(w)^2 r_\lambda^2 + (st)^2 + (s+t)^2 r_\lambda \in \boldsymbol{F}^{\square}(r_\lambda)$, as claimed in assertion b.

We now use the orthogonal basis $v_1, v_2, v_3, v_4$, again. If the discriminant of $h$ is trivial then $abc$ is a square in $\boldsymbol{F}$, and with respect to the basis $v_1, v_2, v_3, \sqrt{abc}\, c^{-1} v_4$, the form $h$ has the diagonal Gram matrix with diagonal entries $1, a, b, ab$. The endomorphisms $A$ and $B$ from above then give similitudes of $h$, with multipliers $a$ and $b$, respectively. This yields $\boldsymbol{L}_h = \boldsymbol{L}$ and $r(\boldsymbol{L}_h) = r(\boldsymbol{L}) = q(V)$, so assertion c is proved.

Now assume that $q(V)$ is a subfield of $\boldsymbol{F}$, so $c \in \boldsymbol{F}^{\square}(a, b)$. The discriminant of $h$ is represented by $abc$, and by assertion b that element lies in each subfield $\boldsymbol{F}^{\square}(r_\delta)$, where $\delta \in \mathrm{GO}(h) \smallsetminus \boldsymbol{F}\boldsymbol{1}$. If the discriminant is not trivial then $abc$ is not a square, and we have $\{r_\delta \mid \delta \in \mathrm{GO}(h)\} \subseteq \boldsymbol{F}^{\square}(abc)$. Again, there exist $s_1, s_a, s_b, u \in \boldsymbol{F}$ with $c = s_1^2 + s_a^2 a + s_b^2 b + u^2 ab$, and $u \neq 0$ because $1, a, b, c$ are linearly independent over $\boldsymbol{F}^{\square}$. Replacing the last basis vector $v_4$ by $u^{-1} v_4$, we achieve $u = 1$, and

$$(*) \qquad ab = s_1^2 + s_a^2 a + s_b^2 b + c\,.$$

We claim that there exists $\gamma \in \mathrm{GO}(V, h)$ with multiplier $r_\gamma = abc$. In order to find the images of the basis vectors $v_1, v_2, v_3, v_4$ under $\gamma$, we use the relation $(*)$ to compute

$$\begin{aligned} q(\gamma(v_1)) &= abc &&= (ab + s_1^2)^2 + (s_b b + s_1 s_a)^2 a + (s_a a + s_1 s_b)^2 b + s_1^2 c\,, \\ q(\gamma(v_2)) &= a^2 bc &&= a^2\big((s_b b + s_1 s_a)^2 + (b + s_a^2)^2 a + (s_1 + s_a s_b)^2 b + s_a^2 c\big)\,, \\ q(\gamma(v_3)) &= ab^2 c &&= b^2\big((s_a a + s_1 s_b)^2 + (s_1 + s_a s_b)^2 a + (a + s_b^2)^2 b + s_b^2 c\big)\,, \\ q(\gamma(v_4)) &= abc^2 &&= c^2\big(s_1^2 + s_a^2 a + s_b^2 b + c\big)\,. \end{aligned}$$

As $q$ is injective, this yields

$$\begin{aligned}
\gamma(v_1) &= (ab+s_1^2)v_1+(s_bb+s_1s_a)v_2+(s_aa+s_1s_b)v_3+s_1v_4\,,\\
\gamma(v_2) &= a\big((s_bb+s_1s_a)v_1+(b+s_a^2)v_2+(s_1+s_as_b)v_3+s_av_4\big)\,,\\
\gamma(v_3) &= b\big((s_aa+s_1s_b)v_1+(s_1+s_as_b)v_2+(a+s_b^2)v_3+s_bv_4\big)\,,\\
\gamma(v_4) &= c\big(s_1v_1+s_av_2+s_bv_3+v_4\big)\,.
\end{aligned}$$

Straightforward computations (using the relation $(*)$, again) yield that $\gamma(v_1),\gamma(v_2),\gamma(v_3),\gamma(v_4)$ is an orthogonal basis with respect to $h$, and thus $\gamma$ is indeed a similitude with multiplier $abc$. We obtain $\boldsymbol{L}_h=\boldsymbol{F}\mathbf{1}+\boldsymbol{F}\gamma$ and the proof of assertion d is complete.

If $q(V)$ is not a subfield of $\boldsymbol{F}$ then assertion a yields $\boldsymbol{L}=\boldsymbol{F}\mathbf{1}$, and $\boldsymbol{L}_h=\boldsymbol{L}$ follows from $\boldsymbol{F}\mathbf{1}\le\boldsymbol{L}_h\le\boldsymbol{L}$. The discriminant is non-trivial by assertion c, and the last assertion e is established. □

We collect the results obtained for the different cases:

**3.5 Theorem.** *Assume* $\operatorname{char}\boldsymbol{F}=2$ *and* $\ell=2$. *Then one of the following holds.*

- **a.** *If $q$ has defect 3, then $\boldsymbol{K}$ splits and $\mathrm{O}(V,h)\cong\big(\mathrm{SL}_2(\boldsymbol{F})\ltimes\boldsymbol{F}^2\big)\times\boldsymbol{F}$. The normal subgroup $\Xi\cong\boldsymbol{F}^2\times\boldsymbol{F}$ is the kernel of the action on $Wz$. See 3.1.*
  *We obtain $\mathrm{SL}_2(\boldsymbol{F})\cong\eta^o(\mathrm{O}(V,h))<\mathrm{O}(W,g)\cong\mathrm{SL}_2(\boldsymbol{K})\cong\mathrm{SL}_2(\boldsymbol{F})\times\mathrm{SL}_2(\boldsymbol{F})$.*
  *Every multiplier is a square, so $\boldsymbol{L}_h=\boldsymbol{L}=\boldsymbol{F}\mathbf{1}$, and $\mathrm{GO}(V,h)=\boldsymbol{F}^\times\,\mathrm{O}(V,h)$.*
- **b.** *If $q$ has defect 2 and $\boldsymbol{K}$ is not split, then $\mathrm{O}(V,h)\cong\mathrm{SL}_2(\boldsymbol{F})$. See 3.2.a.*
  *We have $\mathrm{SL}_2(\boldsymbol{F})\cong\eta(\mathrm{O}(V,h))<\mathrm{O}(W,g)\cong\mathrm{SL}_2(\boldsymbol{K})$.*
  *The field $\boldsymbol{L}_h=\boldsymbol{L}$ is an inseparable quadratic extension field of $\boldsymbol{F}\mathbf{1}$.*
- **c.** *If $q$ has defect 2 and $\boldsymbol{K}$ is split, then $\mathrm{O}(V,h)\cong\boldsymbol{F}^3$ is abelian, and its action on $Wz$ is neither trivial nor faithful. See 3.2.b.*
  *We have $\mathrm{O}(W,g)\cong\boldsymbol{F}\times\big((\boldsymbol{F}^3\times\boldsymbol{F}^2)\rtimes\boldsymbol{F}\big)$.*
  *Again, the field $\boldsymbol{L}_h=\boldsymbol{L}$ is an inseparable quadratic extension field of $\boldsymbol{F}\mathbf{1}$.*
- **d.** *If $q$ has defect 1, then $\mathrm{O}(V,h)\cong\boldsymbol{F}$ is abelian, and $\boldsymbol{K}$ is not split. See 3.3.*
  *We have $\mathrm{O}(V,h)\cong(\boldsymbol{F},+)$, and $\eta(\mathrm{O}(V,h))<\mathrm{O}(W,g)\cong(\boldsymbol{K},+)$.*
  *The field $\boldsymbol{L}_h$ equals $\boldsymbol{F}\mathbf{1}$, and the group of similitudes is $\mathrm{GO}(V,h)=\boldsymbol{F}^\times\,\mathrm{O}(V,h)$.*
- **e.** *If $q$ has defect 0, then $q$ is anisotropic and both $\mathrm{O}(V,q)$ and $\mathrm{O}(V,h)$ are trivial. See 3.4.*
  - *If $q(V)$ is a subfield of $\boldsymbol{F}$ then $\boldsymbol{L}$ is a totally inseparable extension of degree 4 and exponent 1 over $\boldsymbol{F}\mathbf{1}$, and $r(\boldsymbol{L})=q(V)$.*
  - *If $q(V)$ is not a subfield then the discriminant is not trivial, and $\boldsymbol{L}_h=\boldsymbol{L}=\boldsymbol{F}\mathbf{1}$.*
  - *If the discriminant is trivial then $q(V)=r(\boldsymbol{L})$, and $\boldsymbol{L}_h=\boldsymbol{L}$ has degree 4 over $\boldsymbol{F}\mathbf{1}$.*
  - *If the discriminant is not trivial but $q(V)$ is a field then $\boldsymbol{L}_h$ is an inseparable extension of degree 2 over $\boldsymbol{F}\mathbf{1}$, and $\boldsymbol{F}\mathbf{1}\lneqq\boldsymbol{L}_h\lneqq\boldsymbol{L}$.* □

## Declaration of competing interest

There is no competing interest.

## Acknowledgements

This research was funded by the Deutsche Forschungsgemeinschaft through a Polish-German *Beethoven* grant KR1668/11, and under Germany's Excellence Strategy EXC 2044-390685587, Mathematics Münster: Dynamics-Geometry-Structure.

The authors thank the anonymous referee for helpful suggestions regarding the exposition and for pointing out an error in a previous version of 1.5c. This leads to a substantial improvement of the paper.

## Data availability

No data was used for the research described in the article.